\documentstyle[12pt]{article}
\font\f=cmr10

\begin{document}
$$
$$
\begin{center}
{\bf KAPPA-DEFORMED SPACE-TIME UNCERTAINTY
RELATIONS}\footnote{Supported by KBN grant 2P30208706}

\vspace*{2ex}
{\f ANATOL NOWICKI}

\vspace*{2ex}
{\f Institute of Physics, Pedagogical University,}\\
{\f Plac S{\l}owia\'{n}ski 6, 65-029 Zielona G\'{o}ra, Poland}\\[1cm]
{ January 1997}
\end{center}

\vspace*{4ex}
\begin{center}
\begin{minipage}[]{11.6cm}
{\f ABSTRACT:  We discuss the $\kappa$-deformed phase space obtained
as a cross product algebra of the deformed translations algebra and
its dual configuration space. We consider two kinds of the
$\kappa$-deformed uncertainty relations.}
\end{minipage}
\end{center}

\vspace*{4ex}
{\bf 1. Introduction}

\vspace*{4ex}
Deformations of space-time symmetries are extensively investigated
in last years. In this approach the notion of symmetries is
generalized to quantum groups i.e. Hopf algebras and in consequence
we deal with noncommuting space-time. Having a deformed configuration
space  it is interesting to build and investigate a deformed phase
space. For a given configuration space one can define a momentum
space using the concept of duality. It appears that such a momentum
space has also a Hopf algebra structure.\\
The problem arises when we need to define a deformed phase space as a
generalization of the deformed configuration and momentum spaces.
Roughly speaking, one can define the commutators between position
and momentum in different ways. Therefore, the generalization
procedure is ambiguous.

It is known that from a pair of dual Hopf algebras describing the
quantum symmetry group and its dual quantum Lie algebra one can
construct three different double algebras: Drinfeld double
and Drinfeld codouble, both with Hopf algebra structure, and
Heisenberg  double (or more generally, cross product algebra),
which is not a Hopf algebra.

It is easy to see that deformed phase space with a Drinfeld double
structure does not give us in the nondeformed limit the standard
quantum mechanical phase space, because the canonical commutation
relations between momentum and position operators are broken.
Therefore,  deformed phase space cannot be a Hopf algebra.

In the present note we consider a deformed phase space which is an
extension of the $\kappa$-deformed Minkowski space [1], dual to
the momentum sector of $\kappa$-Poincar\'{e} algebra [2].

Choosing different realizations of $\kappa$-Poincar\'{e} algebra
we can obtain different $\kappa$-deformed phase spaces by the cross
product algebra construction. In particular, we discuss two
realizations of $\kappa$-Poincar\'{e} basis: {\it bicrossproduct
basis} [4] and {\it standard basis} [2].

\vspace*{4ex}
{\bf 2. $\kappa$-deformed phase space as cross product algebra}

\vspace*{4ex}
Let us begin our considerations assuming the bicrossproduct basis
for $\kappa$-Poincar\'{e} algebra [4], then the $\kappa$-deformed
Hopf algebra of translations ${\cal{P}}_{\kappa}$
is given by $(k=1,2,3)$\medskip\\
i) {\it bicrossproduct basis}\medskip

\hspace*{4em}$[P_{0}, P_{k}] = 0\hfill(1.a)$\bigskip

\hspace*{4em}$\Delta(P_0) = P_0 \otimes 1 + 1\otimes P_0$\medskip

\hspace*{4em}$\Delta(P_k) = P_k \otimes 1 + e^{-{{P_{0}}\over
\kappa\hbar}}\otimes P_k\hfill(1.b)$\medskip\\
and we define the antipode and counit as follows \medskip

\hspace*{4em}$S(P_\mu)= -P_\mu\qquad\qquad
\epsilon(P_\mu)=0\hfill(2)$\medskip\\
where two constants are introduced: $\hbar$ - Planck's
constant  and $\kappa$ - massive deformation parameter.

Using the duality relations $(\mu,\nu=0,1,2,3)$
$$
<x_\mu, P_\nu> = i\hbar g_{\mu\nu}\qquad g_{\mu\nu} = (-1,1,1,1)
\eqno(3)
$$
we obtain the noncommutative $\kappa$-deformed configuration space
${\cal{X}}_{\kappa}$ as a Hopf algebra with the following algebra and
coalgebra structure [1]\bigskip

\hspace*{4em}$[x_{0}, x_{k}] = {i\over \kappa}x_k$\medskip

\hspace*{4em}$[x_{k}, x_{l}] = 0\hfill(4a)$\bigskip

\hspace*{4em}$\Delta(x_\mu ) = x_{\mu}\otimes 1 + 1\otimes
x_{\mu}$\medskip

\hspace*{4em}$S(x_{\mu}) = -x_{\mu}\qquad\qquad \epsilon(x_{\mu}) =
0\hfill(4b)$\bigskip\\
It is obvious that the algebra ${\cal{P}}_{\kappa}$ is a left
(right) ${\cal{X}}_{\kappa}$-module and vice versa the algebra
${\cal{X}}_{\kappa}$ is
${\cal{P}}_{\kappa}$-module if we introduce the following actions
(we use the Sweedler's notation)\medskip\\
-{\it left actions}
$$
\triangleright : {\cal{X}}_{\kappa}\otimes
{\cal{P}}_{\kappa}\rightarrow
{\cal{P}}_{\kappa}: x\otimes p\rightarrow x\triangleright
p=<x, p_{(2)}>p_{(1)} \eqno(5.a)
$$
$$
\triangleright : {\cal{P}}_{\kappa}\otimes
{\cal{X}}_{\kappa}\rightarrow
{\cal{X}}_{\kappa}: p\otimes x\rightarrow
p\triangleright x=<p,x_{(2)}>x_{(1)} \eqno(5.b)
$$
-{\it right actions}
$$
\triangleleft : {\cal{P}}_{\kappa}\otimes
{\cal{X}}_{\kappa}\rightarrow
{\cal{P}}_{\kappa}: p\otimes x\rightarrow p\triangleleft
x=<x, p_{(1)}>p_{(2)} \eqno(5.c)
$$
$$
\triangleleft : {\cal{X}}_{\kappa}\otimes
{\cal{P}}_{\kappa}\rightarrow
{\cal{X}}_{\kappa}: x\otimes
p\rightarrow x\triangleleft  p=<p,x_{(1)}>x_{(2)} \eqno(5.d)
$$
in choosen basis the actions are the following\medskip\\
$\triangleright : {\cal{X}}_{\kappa}\otimes
{\cal{P}}_{\kappa}\rightarrow {\cal{P}}_{\kappa}:$
$$
\begin{array}{ll}
x_{0}\triangleright P_{0} \ = \ -i\hbar\qquad &  x_{k}\triangleright
P_{0} \ = \ 0\medskip\\
x_{0}\triangleright P_{k} \ = \ 0\qquad &  x_{k}\triangleright P_{l}
\ = \ i\hbar\delta_{kl}e^{-{P_{0}\over \kappa\hbar}}
\end{array}
\eqno(6.a)
$$
$\triangleright : {\cal{P}}_{\kappa}\otimes
{\cal{X}}_{\kappa}\rightarrow {\cal{X}}_{\kappa}:$
$$
\begin{array}{ll}
P_{0}\triangleright x_{0} \ = \ -i\hbar\qquad &  P_{k}\triangleright
x_{0} \ = \ 0\medskip\\
P_{0}\triangleright x_{k} \ = \ 0\qquad &  P_{k}\triangleright
x_{l} \ = \ i\hbar\delta_{kl}
\end{array}
\eqno(6.b)
$$
$\triangleleft : {\cal{P}}_{\kappa}\otimes {\cal{X}}_
{\kappa}\rightarrow {\cal{P}}_{\kappa}:$
$$
\begin{array}{ll}
P_{0}\triangleleft x_{0} \ = \ -i\hbar\qquad &  P_{k}\triangleleft
x_{0} \ = \ {i\over \kappa}P_{k}\medskip\\
P_{0}\triangleleft x_{k} \ = \ 0\qquad &  P_{k}\triangleleft x_{l}
\ = \ i\hbar\delta_{kl}
\end{array}
\eqno(6.c)
$$
$\triangleleft : {\cal{X}}_{\kappa}\otimes
{\cal{P}}_{\kappa}\rightarrow {\cal{X}}_{\kappa}:$
$$
\begin{array}{ll}
x_{0}\triangleleft P_{0} \ = \ -i\hbar\qquad &  x_{k}\triangleleft
P_{0} \ = \ 0\medskip\\
x_{0}\triangleleft P_{k} \ = \ 0\qquad &  x_{k}\triangleleft P_{l}
\ = \ i\hbar\delta_{kl}
\end{array}
\eqno(6.d)
$$

In order to construct a deformed phase space $\Pi_{\kappa}$ one
has to extend the commutation relations $(1a)$
and $(4b)$ by adding a cross
commutators between ${\cal{X}}_{\kappa}$ and ${\cal{P}}_{\kappa}$.
In this way  $\Pi_{\kappa}
\sim {\cal{X}}_{\kappa}\otimes {\cal{P}}_{\kappa}$ as a vector space
and becomes an associative algebra.\\
We find the cross commutators using the notion of a left (right)
{\it cross product (smash product) algebra}.\medskip\\
We recall the following  definition\medskip\\
{\bf Def.}({\it cross product algebra}) [5]\medskip

Let ${\cal{P}}$ be a Hopf algebra and ${\cal{X}}$ a left (right)
${\cal{P}}$-module algebra. Left (right) {\it cross
product algebra} ${\cal{X}} >\!\!\!\triangleleft {\cal{P}}$
(${\cal{P}}\triangleright\!\!\!<{\cal{X}}$)
is a vectors space ${\cal{X}}\otimes {\cal{P}}$ with product
$$
(x\otimes p)(\tilde{x}\otimes \tilde{p})=x(p_{(1)}\triangleright
\tilde{x})\otimes p_{(2)}\tilde{p}\qquad(left\,\,\, cross\,\,\,
product) \eqno(7.a)
$$
and
$$
(p\otimes x)(\tilde{p}\otimes \tilde{x})=p\tilde{p}_{(1)}\otimes
(x\triangleleft \tilde{p}_{(2)})\tilde{x}\qquad(right\,\,\,
cross\,\,\, product)
\eqno(7.b)
$$
with unit element $1\otimes 1$, where $x, \tilde{x}\in {\cal{X}}$
and $p, \tilde{p}\in {\cal{P}}$.\medskip\\
In our case the notion of cross product algebra is equivalent to
the Heisenberg double algebra [3].

The obvious isomorphism ${\cal{X}}\sim {\cal{X}}\otimes 1$,
${\cal{P}}\sim
{\cal{P}}\otimes 1$ gives us the following cross relations
between the configuration and momentum space\medskip

\hspace*{4em}$p\circ x=(p_{(1)}\triangleright x)\circ p_{(2)}$
\qquad\qquad {\it (left cross product)}\hfill$(8.a)$\medskip

\hspace*{4em}$x\circ p=p_{(1)}\circ (x\triangleleft p_{(2)})$
\qquad\qquad {\it (right cross product)}\hfill$(8.b)$\medskip\\
which one can rewrite as a cross commutation relations  $[x,p]=x\circ
p-p\circ x$.\medskip

Let us notice that the following isomorphisms hold between left and
right cross product algebras\medskip

\hspace*{4em}${\cal{X}}_{\kappa}>\!\!\!\triangleleft
{\cal{P}}_{\kappa}\sim
{\cal{X}}_{\kappa}\triangleright\!\!\!<{\cal{P}}_{\kappa}
\hfill(9.a)$\medskip

\hspace*{4em}${\cal{P}}_{\kappa}>\!\!\!\triangleleft
{\cal{X}}_{\kappa}\sim
{\cal{P}}_{\kappa}\triangleright\!\!\!<{\cal{X}}_{\kappa}\hfill(9.b)
$\medskip\\ therefore, it is enough to consider left cross product
algebras ${\cal{X}}_{\kappa}>\!\!\!\!\triangleleft {\cal{P}}_{\kappa}$
and ${\cal{P}}_{\kappa}>\!\!\!\triangleleft {\cal{X}}_{\kappa}$ as
the models of our $\kappa$-deformed phase space $\Pi_{\kappa}$.

From the relations $(1.a)$, $(4.a)$ and the cross relations
$(7.a)$ it follows that we have the commutation relations in the
form\medskip\\ - {\it the case} $\Pi_{\kappa}\sim {\cal{X}}_
{\kappa}>\!\!\!\triangleleft {\cal{P}}_{\kappa}$\bigskip

\hspace*{4em}$[P_{0}, P_{k}] \ = \ 0\hfill(10.a)$\bigskip

\hspace*{4em}$[x_{\mu}, x_{\nu}] \ = \ {i\over
\kappa}(\delta_{\mu0}x_{\nu} - \delta_{\nu0}x_{\mu})\hfill(10.b)
$\bigskip

\hspace*{4em}$[x_{\mu}, P_{0}] \ = \ i\hbar \delta_{\mu
0}\hfill(10.c)$\bigskip

\hspace*{4em}$[x_{\mu}, P_{k}] \ = \ -i\hbar \delta_{\mu k} \ -
\ {i\over \kappa}\delta_{\mu 0} \ P_{k}\hfill(10.d)$\bigskip

In analogous way, from $(1.a)$, $(4.a)$ and $(7.b)$ we
obtain\medskip\\ - {\it the case} $\Pi_{\kappa}\sim
{\cal{P}}_{\kappa}>\!\!\!\triangleleft {\cal{X}}_{\kappa}$\bigskip

\hspace*{4em}$[P_{0}, P_{k}] \ = \ 0\hfill(11.a)$\bigskip

\hspace*{4em}$[x_{\mu}, x_{\nu}] \ = \ {i\over
\kappa}(\delta_{\mu0}x_{\nu} - \delta_{\nu0}x_{\mu})\hfill(11.b)
$\bigskip

\hspace*{4em}$[x_{0}, P_{\mu}] \ = \ -i\hbar \delta_{\mu0}
\hfill(11.c)$\bigskip

\hspace*{4em}$[x_{k}, P_{\mu}] \ = \ i\hbar \delta_{\mu k} \
e^{-{P_{0}\over \kappa\hbar}}\hfill(11.d)$\bigskip\\
If we choose $g_{\mu\nu}$ with oposite signs in $(3)$ then the
relations $(10.b-d)$ and $(11.b-d)$ change signs.\\
On the other hand, if we choose instead $(1.b)$ transposed coproduct
$$
\tilde{\Delta}(P_k) \ = \ P_{k}\otimes e^{-{{P_0}\over \kappa\hbar}}
+ 1\otimes P_{k}
$$
the commutators $(10)$ change signs and we obtain $\kappa$-deformed
phase space discussed in [3].

For our choice we obtain in the limit $\kappa\rightarrow\infty$ the
quantum-mechanical phase space in the case
$\Pi_{\kappa}\sim {\cal{P}}_{\kappa}>\!\!\!\triangleleft
{\cal{X}}_{\kappa}$ (relations $(11)$).\bigskip

If we consider the $\kappa$-Poincar\'{e} algebra in the standard
basis [2], the translations subalgebra  $\Pi_\kappa$ is given
by\medskip\\ ii) {\it standard basis}\medskip

\hspace*{4em}$[P_{0}, P_{k}] = 0\hfill(12.a)$\bigskip

\hspace*{4em}$\Delta(P_0) = P_0 \otimes 1 + 1\otimes P_0$\medskip

\hspace*{4em}$\Delta(P_k) = P_{k}\otimes e^{{P_{0}\over
2\kappa\hbar}} + e^{-{P_{0}\over 2\kappa\hbar}}\otimes
P_{k}\hfill(12.b)$\medskip

In analogous way we obtain two cross product algebras\medskip\\
- {\it the case} $\Pi_{\kappa}\sim {\cal{X}}_{\kappa}>
\!\!\!\triangleleft {\cal{P}}_{\kappa}$\bigskip

\hspace*{4em}$[P_{0}, P_{k}] \ = \ 0\hfill(13.a)$\bigskip

\hspace*{4em}$[x_{\mu}, x_{\nu}] \ = \ {i\over
\kappa}(\delta_{\mu0}x_{\nu} - \delta_{\nu0}x_{\mu})
\hfill(13.b)$\bigskip

\hspace*{4em}$[x_{\mu}, P_{0}] \ = \ i\hbar \delta_{\mu 0}
\hfill(13.c) $\bigskip

\hspace*{4em}$[x_{\mu}, P_{k}] \ = \ -i\hbar \delta_{k \mu}
e^{-{P_{0}\over 2\kappa\hbar}} \ - \ {i\over 2\kappa}
\delta_{0 \mu}P_{k} \hfill(13.d)$\bigskip\\
- {\it the case} $\Pi_{\kappa}\sim {\cal{P}}_{\kappa}>
\!\!\!\triangleleft {\cal{X}}_{\kappa}$\bigskip

\hspace*{4em}$[P_{0}, P_{k}] \ = \ 0\hfill(14.a)$\bigskip

\hspace*{4em}$[x_{\mu}, x_{\nu}] \ = \ {i\over
\kappa}(\delta_{\mu0}x_{\nu} - \delta_{\nu0}x_{\mu})\hfill(14.b)
$\bigskip

\hspace*{4em}$[x_{\mu}, P_{0}] \ = \ -i\hbar \delta_{\mu 0}
\hfill(14.c)$\bigskip

\hspace*{4em}$[x_{\mu}, P_{k}] \ = \ i\hbar \delta_{\mu k} \
e^{-{P_{0}\over
2\kappa\hbar}} - {i\over 2\kappa\hbar}\delta_{\mu 0}P_k\hfill
(14.d)$\bigskip

Let us notice that only for the algebra $(14)$ we get the correct
limit
for $\kappa\rightarrow\infty$, so one can assume that the relations
$(14)$ define $\kappa$-deformed phase space.

\vspace*{4ex}
{\bf 3. $\kappa$-deformed uncertainty relations}

\vspace*{4ex}
Introducing the dispersion of the observable  $a$ in quantum
mechanical sense by\medskip
$$
\Delta(a) \ = \ \sqrt{<a^2> - <a>^2}
\eqno(15.a)
$$
we have
$$
\Delta(a)\Delta(b)\geq {1\over 2}|<c>|\qquad where\qquad c=[a,b]
\eqno(15.b)
$$

In the first case $(11)$, we obtain
$\kappa$-deformed uncertainty relations in $\Pi_{\kappa}$ in the
form\bigskip

\hspace*{4em}$\Delta(x_0)\Delta(x_k)\geq
{1\over2\kappa}|<x_k>|\hfill(16.a)$\bigskip

\hspace*{4em}$\Delta(P_k) \Delta(x_l)\geq {1\over
2}\hbar\delta_{kl}|<e^{-{{P_0}\over k\hbar}}>| \hfill(16.b)$\bigskip

\hspace*{4em}$\Delta(P_0) \Delta(x_0)\geq {1\over
2}\hbar\hfill(16.c)$\bigskip

\hspace*{4em}$\Delta(P_k) \Delta(x_0)\geq 0\hfill(16.d)$\bigskip\\
and in the second case (14) we get\bigskip

\hspace*{4em}$\Delta(x_0)\Delta(x_k)\geq {1\over2\kappa}|<x_k>|
\hfill(17.a)$\bigskip

\hspace*{4em}$\Delta(P_k) \Delta(x_l)\geq {1\over
2}\hbar\delta_{kl}|<e^{-{P_{0} \over 2\kappa\hbar}}>|
\hfill(17.b)$\bigskip

\hspace*{4em}$\Delta(P_0) \Delta(x_0)\geq {1\over 2}\hbar\hfill
(17.c)$\bigskip

\hspace*{4em}$\Delta(P_k) \Delta(x_0)\geq  {1\over 4\kappa}|<P_{k}>|
\hfill(17.d)$\bigskip\\
In both cases we obtain the standard uncertainty relations in the
limit $\kappa \rightarrow \infty$. The physical implications of these
deformed uncertainty relations will be considered elsewhere.

\vspace*{4ex}
{\bf 4. Final remarks}

\vspace*{4ex}
In this note we presented two realizations of $\kappa$-deformed phase
space. The cross relations between position and momentum operators
form two kinds of generalized Heisenberg commutation relations
depending on two parameters: $\hbar$ - Planck's constant
and $\kappa$ - deformation parameter with mass dimension.

We would like to stress that the realizations of the cross product
algebras of $\kappa$-deformed translations ${\cal{P}}_{\kappa}$ and
positions ${\cal{X}}_{\kappa}$ strongly depend on the choice of the
basis of $\kappa$-Poincar\'{e} algebra as well on the metric tensor
$g_{\mu\nu}$ (see also [3]).

\vspace*{4ex}
{\bf Acknowledgments}

\vspace*{4ex}
I would like to thank J. Lukierski for interesting discussions on
this subject and critical remarks.

\vspace*{4ex}
{\bf Bibliography}

\vspace*{4ex}
\itemsep=-.2pc
\begin{enumerate}
\item S.~Zakrzewski, {\em Journ. Phys.} {\bf A 27} (1994) 2075.\\
P.~Kosinski, P.~Maslanka, ``The Duality Between
$\kappa$-Poincar\'{e} Algebra and $\kappa$-Poincar\'{e} Group'',
Lodz preprint IM UL 3/1994.
\item J.~Lukierski, A.~Nowicki, H.~Ruegg, {\em Phys. Lett.} {\bf B
293} (1993) 419.
\item J.~Lukierski, A.~Nowicki, ``Heisenberg Double Description of
$\kappa$-Poincar\'{e} Algebra and $\kappa$-deformed Phase Space'' in:
{\it Proceedings of XXI International Colloquium on Group Theoretical
Methods in Physics}, Goslar 1996.
\item S.~Majid, H.~Ruegg, {\em Phys. Lett.} {\bf B 334} (1994) 348.
\item S.~Majid, ``{\it Foundations of Quantum Group Theory}'',
Cambridge University Press, Cambridge 1995.
\item J.~Lukierski, H.~Ruegg, W.J.~Zakrzewski, {\em Ann. of Phys.}
{\bf 243} (1995) 90.
\end{enumerate}

\end{document}